\documentclass[letterpaper, 10 pt, conference]{ieeeconf}  

\IEEEoverridecommandlockouts                              
\usepackage{hyperref}
\hypersetup{%
  colorlinks=true, 	
  linktoc=all,     	
  linkcolor=black,	
  citecolor= black,
  urlcolor=blue
}

\usepackage{epsfig} 
\usepackage[table]{xcolor}
\usepackage{amsmath} 
\usepackage{amssymb}  
\usepackage{subcaption}
\usepackage{mcode}
\usepackage{url}
\usepackage{listings}
\usepackage{balance}
\usepackage{mathtools}
\newcommand{\reals}{\mathbf{R}}

\newcommand{\ones}{\mathbf{1}}

\newcommand{\ie}{\emph{i.e.}}
\newcommand{\eg}{\emph{e.g.}}
\newcommand{\norm}[1]{\lVert#1\rVert}

\title{\LARGE \bf
Embedded Code Generation with CVXPY
}

\author{Maximilian Schaller, Goran Banjac, Steven Diamond, Akshay Agrawal, Bartolomeo Stellato, and Stephen Boyd
\thanks{This work was supported by the European Research Council (ERC) under the European Union's Horizon 2020 research and innovation programme grant agreement OCAL, No. 787845, Stanford's SystemX, and the AI Chip Center for Emerging Smart 
Systems (ACCESS).}
\thanks{M.\ Schaller and G.\ Banjac are with the Department of Information Technology and Electrical Engineering, ETH Zürich, 8092 Zürich, Switzerland.
	{\tt\small\{mschaller, gbanjac\}@ethz.ch}}
\thanks{S.\ Diamond is with Gridmatic, Campbell CA 95008, USA.
	{\tt\small steven@gridmatic.com}}
\thanks{A.\ Agrawal and S.\ Boyd are with the Department of Electrical Engineering, Stanford University, Stanford CA 94305, USA.
	{\tt\small\{akshayka, boyd\}@stanford.edu}}
\thanks{B.\ Stellato is with the Department of Operations Research and Financial Engineering, Princeton University, Princeton NJ 08544, USA.
	{\tt\small bstellato@princeton.edu}}
}

\begin{document}

\maketitle
\thispagestyle{empty}
\pagestyle{empty}

\begin{abstract}
We introduce CVXPYgen, a tool for generating custom C code, suitable for embedded applications,
that solves a parametrized class of convex optimization problems.
CVXPYgen is based on CVXPY, a Python-embedded domain-specific language that supports a natural
syntax (that follows the mathematical description) for specifying convex optimization problems.
Along with the C implementation of a custom solver, CVXPYgen creates a Python wrapper
for prototyping and desktop (non-embedded) applications.
We give two examples, position control of a quadcopter and back-testing a portfolio
optimization model.
CVXPYgen outperforms a state-of-the-art code generation tool in terms of
problem size it can handle, binary code size, and solve times.
CVXPYgen and the generated solvers are open-source.
\end{abstract}

\section{Introduction}

Convex optimization is used in many domains, including signal and image processing~\cite{mattingley2010real,zibulevsky2010l1}, control~\cite{wang2009fast,garcia1989model}, and finance~\cite{boyd2017multi,markowitz1952portfolio}, to mention just a few.
A (parametrized) convex optimization problem can be written as
\begin{equation}\label{eq:cvx}
	\begin{array}{ll}
		\text{minimize} \quad & f_0(x, \theta) \\
		\text{subject to} \quad & f_i(x, \theta) \le 0, \quad i=1,\ldots,p \\
		& g_j(x,\theta) = 0, \quad j=1,\ldots, r,
	\end{array}
\end{equation}
where $x \in \reals^n$ is the optimization variable,
$f_0$ is the objective function to be minimized, $f_1,\ldots,f_p$
are the inequality constraint functions, and $g_1\ldots, g_r$ are the
equality constraint functions.
We require that $f_0, \ldots, f_p$ are convex functions, and
$g_1, \ldots, g_r$ are affine functions~\cite{boyd2004convex}.
The parameter $\theta\in \reals^d$ specifies data that can change, but is constant and given when we solve an instance of the problem.
We refer to the parametrized problem~\eqref{eq:cvx} as
a \emph{problem family}; when we specify a fixed value of $\theta$, we refer to it as a \emph{problem instance}.
We let $x^\star$ denote an optimal point for the problem~\eqref{eq:cvx}, assuming it exists.

The problem family can be specified using a domain-specific language (DSL) for convex optimization.
Such systems allow the user to specify the functions $f_i$ and $g_j$ in a simple
format that closely follows the mathematical description of the problem.
Examples include YALMIP \cite{lofberg2004yalmip} and CVX \cite{grant2014cvx} (in Matlab), CVXPY \cite{diamond2016cvxpy} (in Python),
Convex.jl~\cite{convexjl} and JuMP \cite{Dunning2017jump} (in Julia),
and CVXR~\cite{fu2017cvxr} (in R).
We focus on CVXPY, which also supports the declaration of parameters,
enabling it to specify problem families, not just problem instances.

DSLs parse the problem description and translate (or canonicalize)
it to an equivalent problem
that is suitable for a solver that handles some generic class of problems,
such as linear programs (LPs), quadratic programs (QPs), second-order cone programs (SOCPs), semidefinite programs (SDPs), and others such as exponential cone programs \cite{boyd2004convex}.
We focus on solvers that are suitable for embedded applications, \ie, are
single-threaded, can be statically compiled, and do not make system calls:
OSQP \cite{stellato2020osqp} handles QPs, SCS \cite{odonoghue2016scs} and
ECOS \cite{domahidi2013ecos} handle cone programs that include SOCPs.
After the canonicalized problem is solved, a solution of the original problem is retrieved from a solution of the canonicalized problem.

It is useful to think of the whole process as a function that maps~$\theta$,
the parameter that specifies the problem instance, into $x^\star$, an optimal
value of the variable. With a DSL, this process consists of three steps.
First the original problem description is canonicalized to a problem in some 
standard (or canonical) form;
then the canonicalized problem is solved using a solver; and finally, a solution
of the original problem is retrieved from a solution of the canonicalized problem.

\begin{figure}
	\centering
	\begin{subfigure}{\columnwidth}
		\includegraphics[width=\linewidth]{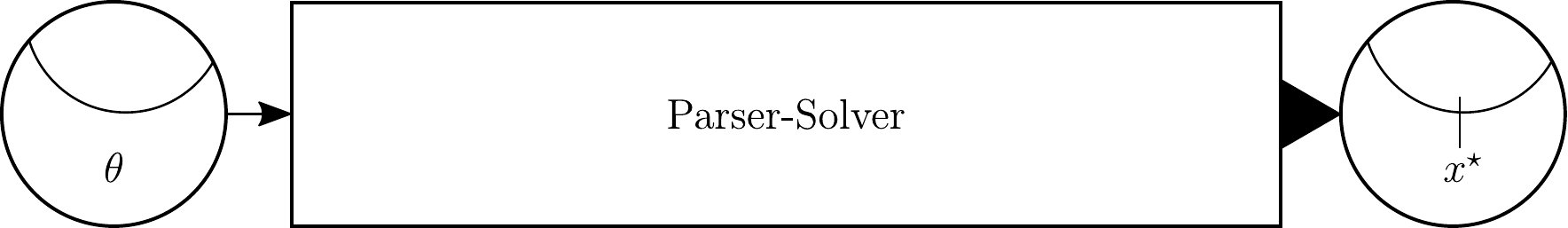}
		\caption{Parser-solver calculating solution $x^\star$ for 
problem instance with parameter $\theta$.}
		\label{fig:flow_general}
	\end{subfigure}
	\par\vspace{10pt}
	\begin{subfigure}{\columnwidth}
		\includegraphics[width=\linewidth]{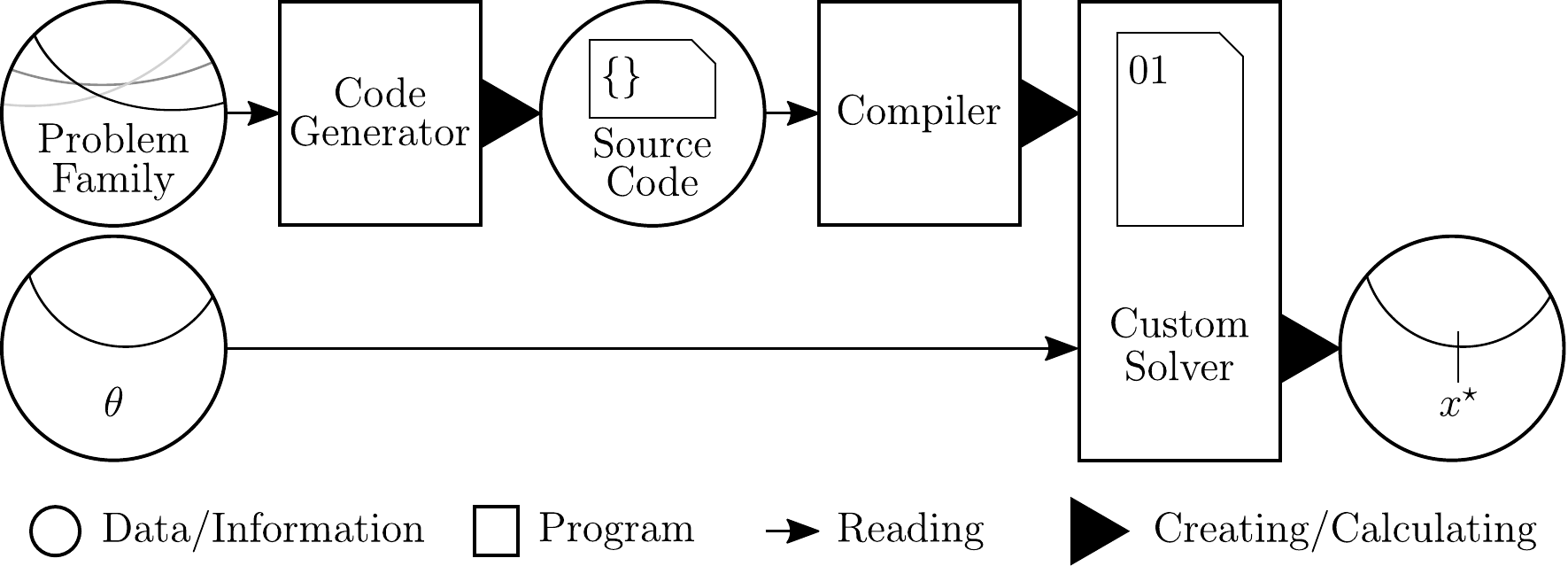}
		\caption{Source code generation for problem family, followed by compilation to 
custom solver. The compiled solver computes a solution $x^\star$ to the problem instance with
parameter $\theta$.}
		\label{fig:flow_custom}
	\end{subfigure}
	\caption{Comparison of convex optimization problem parsing and solving approaches.}
	\vspace{-1em}
	\label{fig:flow}
\end{figure}

Most of these DSLs are organized as
\emph{parser-solvers}, which carry out the canonicalization each time the problem is
solved (with different parameter values).
This simple setting is illustrated in figure~\ref{fig:flow_general}.
We are interested in applications where we solve many instances of the
problem, possibly in an embedded application with hard real-time constraints.
For such applications, a \emph{code generator} makes more sense.
A code generator takes as input a description of a problem family, and generates
specialized solver source code for that specific family.
That source code is then compiled, and we have an efficient solver for
the specific family.

This workflow is illustrated in figure~\ref{fig:flow_custom}.
The compiled solver has a number of advantages over parser-solvers.
First, the compiled solver can be deployed in embedded systems, fulfilling rules for safety-critical code \cite{holzmann2006power}.
Second, by caching canonicalization and exploiting the problem structure, the compiled solver is faster.

A well known code generator is CVXGEN~\cite{mattingley2012cvxgen}.
It handles problems that can be transformed to QPs and includes a custom interior-point solver.
CVXGEN is used in many applications, including autonomous driving, dynamic energy management,
and real-time trading in finance.
SpaceX uses CVXGEN to generate flight code for high-speed onboard convex
optimization for precision landing of space rockets~\cite{blackmore2016autonomous}.

CVXGEN was designed for use in real-time control systems, where the problems solved are
not too big, either in terms of the number of variables or number of parameters.
CVXGEN unrolls for-loops in its generated source code files to increase the solving speed,
but this can also result in large compiled code size.
Due to the flat and explicit code generated, CVXGEN only handles problems with
up to around a few thousand parameters.
(More accurately, CVXGEN is limited to 4000 nonzero entries in the linear system of
equations solved in each iteration.)

\subsection{Contribution}\label{sec:contribution}

In this paper, we introduce the code generation tool CVXPYgen, which produces custom C
code to solve a parametrized family of convex optimization problems.
The design decisions for CVXPYgen are somewhat different from those made for CVXGEN.
First, CVXPYgen is built on top of the DSL CVXPY, whereas CVXGEN is entirely
self-contained.
This means that prototypes can be developed, prototyped, and simulated in Python
using CVXPY.
Second, CVXPYgen interfaces with multiple solvers, currently OSQP, SCS, and ECOS.
This means that CVXPYgen supports problems more general than those that can be
transformed to QPs.
As far as we know, CVXPYgen is the first generic code generator for convex optimization
that supports SOCPs.
When using OSQP or SCS (both based on first-order methods), the generated solvers
support warm-starting, which can bring more speed in some applications~\cite{garcia1989model}.
Third, CVXPYgen does not aggressively unroll loops in the generated code, which allows
it to support high-dimensional parameters.
In addition, matrix parameters can have any user-defined sparsity pattern.
CVXPYgen uses partial update canonicalization, in which only the parameters changed are
processed when solving a new problem instance.
Fourth, CVXPYgen and its generated solvers are fully open-source, whereas CVXGEN is 
proprietary.

CVXPYgen (and more generally, code generation)
is useful for two families of practical applications.
The first is solving convex optimization problems in real-time settings on embedded devices,
as is done in control systems, real-time resource allocators, and other applications.
The second is in solving a large number of instances of a problem family,
possibly on general-purpose computers. One example is back-testing in finance,
where a trading policy based on convex optimization is simulated
on historical or simulated data over many periods. Typical back-tests involve solving
thousands or more instances of a problem family.
In these applications, there is no hard real-time constraint; the goal is simply to
speed up solving by avoiding repeatedly canonicalizing the problem.

\subsection{Prior work}\label{sec:related}
Several other code generators for optimization have been developed in addition to CVXGEN.
FORCESPRO~\cite{FORCESPro} and FORCES NLP~\cite{FORCESNLP} are proprietary code generators for multi-stage control problems.
They handle problems that can be transformed to multi-stage quadratically constrained quadratic programs and
nonlinear programs, respectively.
The open-source code generators QCML~\cite{QCML} and CVXPY-CODEGEN~\cite{cvxpy_codegen},
which interface with ECOS, were
developed before CVXPY included support for parameters.
These early prototypes are no longer actively supported or maintained.

\subsection{Outline}

The remainder of this paper is structured as follows.
In \S\ref{sec:CVXPYgen} we describe, at a high level, how CVXPYgen works, and in
\S\ref{sec:simple_example}, we illustrate how it is used with a simple example.
In \S\ref{sec:comparison_CVXGEN} and \S\ref{sec:comparison_CVXPY} we
compare CVXPYgen to CVXGEN (for embedded use) and CVXPY (for general purpose use), respectively.
We conclude the paper in \S\ref{sec:conclusion}.

\section{CVXPYgen}\label{sec:CVXPYgen}

CVXPYgen is based on the open-source Python-embedded DSL CVXPY.
CVXPY handles  many types of conic programs and certain types of nonconvex problems,
whereas we focus on LPs, QPs, and SOCPs for code generation.
CVXPY provides modeling instructions that follow the mathematical description for convex optimization problems.
It ensures that the modeled problems are convex, using disciplined convex programming (DCP).
DCP is the process of constructing convex functions by assembling given base
functions in mathematical expressions using a simple set of rules \cite{grant2004disciplined}.
DCP ensures that the resulting problem is convex, and also, readily canonicalized
to a standard form.

In DCP, parameters are treated as constants, optionally with specified sign, and there are no
restrictions about how these constants appear in the expressions defining the problem family.
The recently developed concept of disciplined parametrized programming (DPP) puts additional
restrictions on how parameters can enter
a problem description. If a problem family description is DPP-compliant,
then canonicalization and retrieval can be represented as \emph{affine} mappings
\cite{agrawal2019differentiable}.
Thus DPP-compliant problems are reducible to ASA-form, which stands for
\emph{Affine-Solve-Affine} \cite{agrawal2019differentiable}.
This is the key property we exploit in CVXPYgen.
More about the DCP and DPP rules can be found in the aforementioned papers, or at
\url{https://www.cvxpy.org}.

After CVXPY has reduced the DPP-compliant problem to ASA-form, CVXPYgen extracts a
sparse matrix $C$ that canonicalizes the user-defined parameters $\theta$ to
the parameters $\tilde \theta$ appearing in the standard form solver:
\[	\tilde \theta = C \begin{bmatrix}\theta \\ 1 \end{bmatrix}. \]
CVXPYgen analyzes $C$ to determine the user-defined parameters (\ie, components of $\theta$)
that every standardized form parameter depends on.
This information is used when generating the custom solver, where only slices of the
above mapping are computed if not all user-defined parameters are updated between solves.
In addition, it is very useful to know the set of updated canonical parameters when
using the OSQP solver or SCS, as detailed below.

In a similar way the retrieval of the solution $x^\star$ for the original problem from a
solution $\tilde x^\star$ of the canonicalized problem is an affine mapping,
\[
x^\star = R \begin{bmatrix}\tilde x^\star \\ 1 \end{bmatrix},
\]
where $R$ is a sparse matrix.
Typically $R$ is a selector matrix, with only one nonzero entry in each row, equal to one,
in which case this step can be handled via simple pointers in C.

CVXPYgen generates allocation-, library-, and division-free C code for the
canonicalization and retrieval steps, which in essence are nothing more than
sparse matrix-vector multiplication, with some logic that exploits pointers or
partial updates.
Sparse matrices are stored in compressed sparse column format \cite{bulucc2009parallel} and
dense matrices are stored as vectors via column-major flattening.

Any solver can be used to solve the canonicalized problem, which provides the final
link:
\[
\tilde x^\star = \mathcal{S} (\tilde \theta),
\]
where $\mathcal S$ denotes the mapping from the canonicalized parameters to
a solution of the canonicalized problem.
(We assume here that the problem instance is feasible, and that when there are multiple
solutions, we simply pick one.)
If available, CVXPYgen uses the canonical solver's code generation method to
produce C code for canonical solving.
As of now, only OSQP provides this functionality~\cite{banjac2017embedded}.
Otherwise, the solver's C code is simply copied, possibly modified for use
in embedded applications.

OSQP and SCS provide a set of C functions for updating their parameters.
This way, when only canonical vector parameters are updated, the factorization of the linear system involved
in the OSQP or SCS algorithms can be cached and re-used, which can lead to substantial
speed up and division-free code.
In the same way as only the parts of~$\tilde \theta$ are re-canonicalized that depend
on the updated parts of~$\theta$, only the OSQP or SCS update functions associated with
these parameters are called before the canonicalized problem is solved.

The code and the full documentation for CVXPYgen with its generated solvers are available at
\begin{center}
	\url{https://pypi.org/project/cvxpygen}.
\end{center}

\section{Simple example}\label{sec:simple_example}

We consider the nonnegative least squares problem
\begin{equation}\label{eq:nnLS}
	\begin{array}{ll}
		\text{minimize} & \norm{Gx-h}_2^2 \\
		\text{subject to} & x\ge 0,
	\end{array}
\end{equation}
where $x \in \reals^n$ is the variable and
$G \in \reals^{m \times n}$, $h \in \reals^m$ are parameters, so $\theta=(G,h)$.
We will canonicalize this to the standard form accepted by OSQP,
\begin{equation}\label{eq:OSQP}
\begin{array}{ll}
	\text{minimize} \quad & \frac{1}{2} \tilde x^T P \tilde x + q^T \tilde x \\
	\text{subject to} \quad & l \le A \tilde x \le u,
\end{array}
\end{equation}
where $\tilde x \in \reals^{\tilde n}$ is the canonical variable and
all other symbols are canonical parameters, \ie, $\tilde \theta = (P,q,A,l,u)$.
(In this form, entries of $l$ can be $-\infty$, and
entries of $u$ can be $+\infty$.)

The na\"{i}ve canonicalization of \eqref{eq:nnLS} to \eqref{eq:OSQP} takes
$\tilde x=x$ and
\[
P = 2G^TG, \quad q = 2G^Th, \quad A=I, \quad l=0, \quad u=\infty.
\]
In this canonicalization, $\tilde \theta$ is not an affine function of $\theta$,
since some entries of $\tilde \theta$ are products of entries of $\theta$.

The canonicalization that uses DPP first expresses problem~\eqref{eq:nnLS} as
\[
\begin{array}{ll} \text{minimize} & \norm{\tilde x_2}_2^2 \\
\text{subject to} & \tilde x_2 = G\tilde x_1 - h, \quad \tilde x_1\ge 0,
\end{array}
\]
with variable $\tilde x = (\tilde x_1, \tilde x_2)$, where $\tilde x_1=x$
and $\tilde x_2 \in \reals^m$.  We can express this as \eqref{eq:OSQP} with
parameters
\[
P= \begin{bmatrix} 0 & 0 \\ 0 & 2I \end{bmatrix}, \quad
q = 0, \quad
A= \begin{bmatrix} G & -I \\ I & 0 \end{bmatrix},
\]\[
l=(h,0), \quad
u=(h,\infty),
\]
where the second part of $u$ has $\infty$ in every entry, \ie, there is no upper bound on the second part of $A \tilde x$.
In this canonicalization, $\tilde \theta$ is indeed an affine function of
$\theta$.
The retrieval map has the simple (linear) form $x^\star=[I~0] \tilde x^\star$.

\begin{figure}
\lstset{language=Python,
	numbers=left,
	xleftmargin=0.07\columnwidth,
        linewidth=\columnwidth}
\begin{lstlisting}[frame=lines]
import cvxpy as cp
from cvxpygen import cpg

# model problem
x=cp.Variable(n, name='x')
G=cp.Parameter((m,n), name='G')
h=cp.Parameter(m, name='h')
p=cp.Problem(cp.Minimize(cp.sum_squares(G@x-h)),
             [x>=0])
# generate code
cpg.generate_code(p)
\end{lstlisting}
\caption{Code generation for example \eqref{eq:nnLS}. We assume that the dimensions
\texttt{m} and \texttt{n} have been previously defined.}
\label{code:nnLS_CVXPYgen}
\end{figure}
We generate code for this problem as shown in figure~\ref{code:nnLS_CVXPYgen}.
The problem is modeled with CVXPY in lines~5--9.
The actual code generation is done in line~11.
The variable and parameters are named in lines 5--7 via their \texttt{name} attributes.
These names are used for C variable and function naming.

\section{Comparison to CVXGEN}\label{sec:comparison_CVXGEN}

\begin{figure}
	\centering
	\includegraphics[width=\columnwidth]{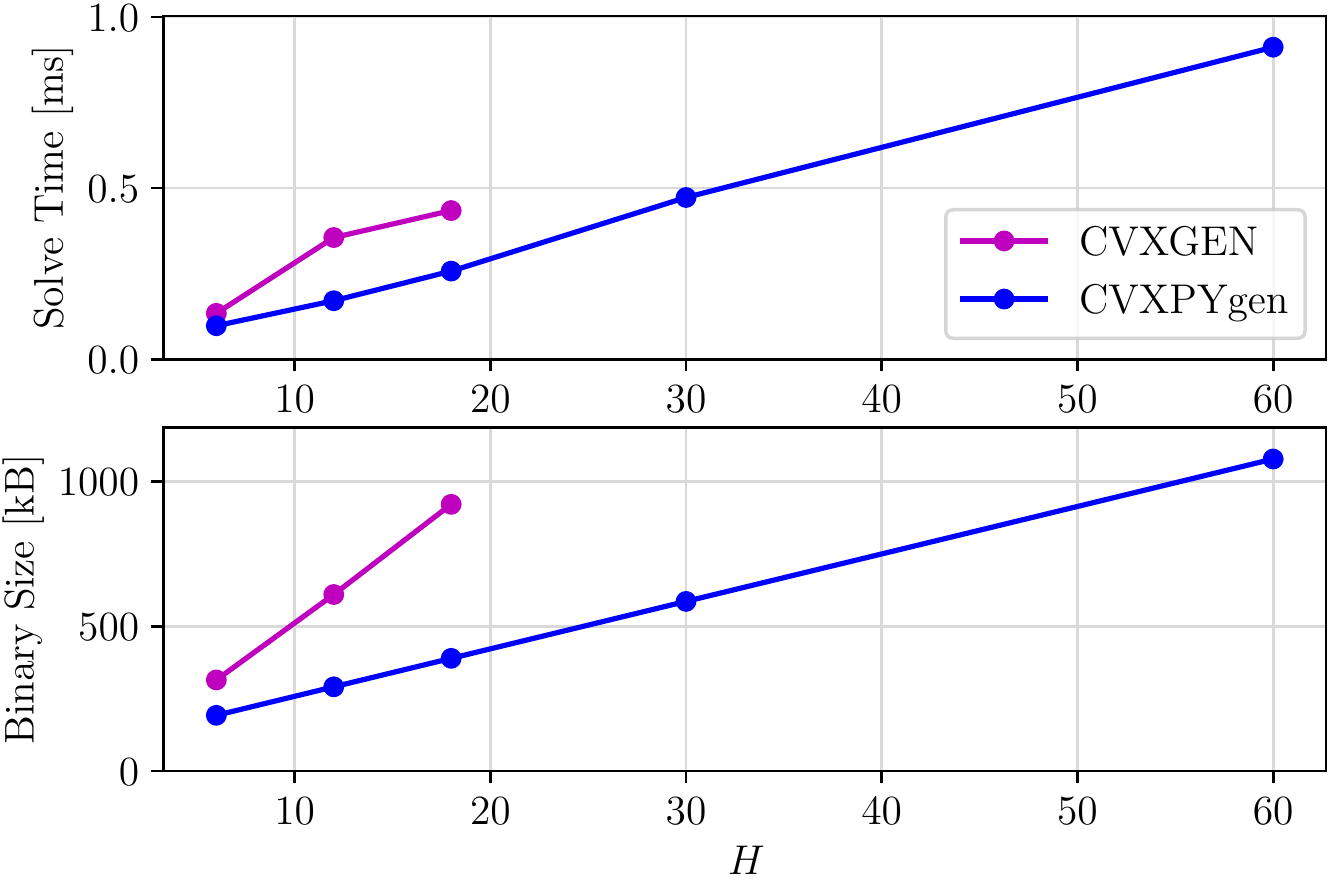}
	\caption{Comparison of solve times (top) and binary sizes (bottom) with CVXGEN (magenta) and CVXPYgen (blue) used for code generation.}
	\label{fig:benchmark_CVXGEN}
\end{figure}

We compare CVXPYgen to CVXGEN for a model predictive control (MPC) problem family,
described in Appendix~\ref{sec:appendix_MPC}.
In particular, we compare solve times of the C interface and executable sizes.
The MPC problems are parametrized by their horizon length
$H\in \left\{6, 12, 18, 30, 60\right\}$; the number of variables is around $10 H$.
Figure~\ref{fig:benchmark_CVXGEN} shows the resulting solve times averaged over
100 simulation steps.
CVXGEN is not able to generate code for $H>18$.
Together with the automatically chosen OSQP solver, CVXPYgen outperforms CVXGEN for all problem sizes.
The bottom of figure~\ref{fig:benchmark_CVXGEN} presents the example executable sizes for CVXGEN and CVXPYgen, respectively.
For all values of $H$, the executables corresponding to CVXPYgen are considerably smaller.

The execution times cited above are on a MacBook Pro 2.3GHz Intel i5.
We have also used these generated solvers to control the position of a custom-built 14-by-14 cm quadcopter.
The generated code was compiled in a robot operating system (ROS) node, and
run on the drone's Intel Atom x5-Z8350 processor, at 30 Hz.
We provide a video of the quadcopter following a circle trajectory at

\begin{center}
\url{https://polybox.ethz.ch/index.php/s/MARR9CGaLqmQaJ0}.
\end{center}

\section{Comparison to CVXPY}\label{sec:comparison_CVXPY}

\begin{figure}
	\centering
	\vspace{4mm}
	\includegraphics[width=\columnwidth]{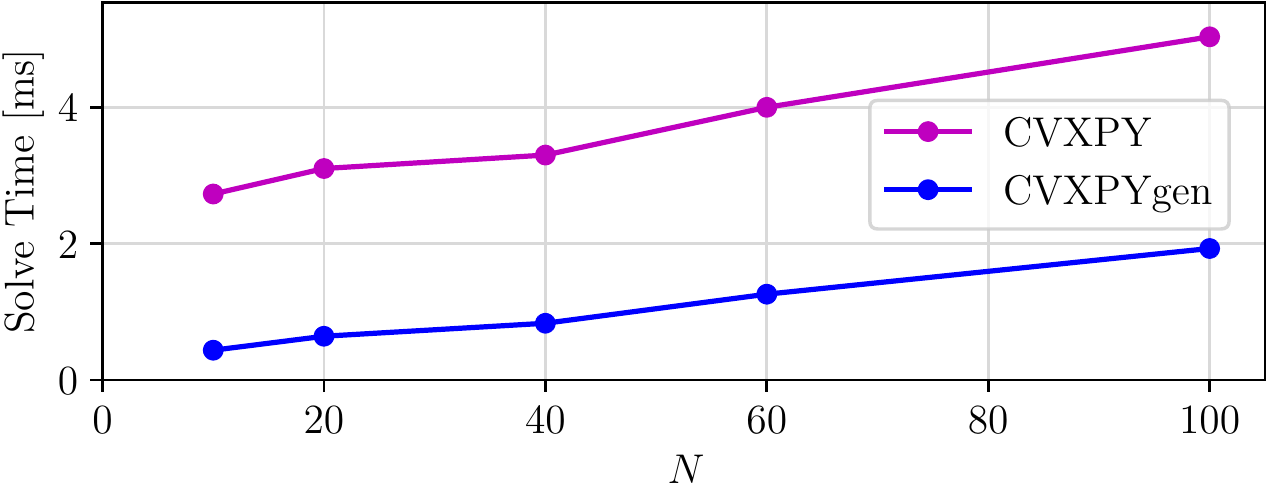}
	\caption{Comparison of solve times with CVXPY (magenta) and the CVXPY
interface of CVXPYgen (blue).}
	\label{fig:benchmark_CVXPY}
\end{figure}

Here we compare CVXPYgen to CVXPY, for a (financial) portfolio optimization problem family,
described in Appendix \ref{sec:appendix_portfolio}.
This family of problems is parametrized by the number of assets in the portfolio
$N\in \{10,20,40,60,100\}$.
The number of variables in these problems is around $2N$.
Figure~\ref{fig:benchmark_CVXPY} gives the results for 500 solves, a two year back-test
using historical data.
We see that the average solver speed is about 6 times faster with CVXPYgen,
for $N=10$, with the ratio dropping to 2.5 for $N=100$.
The execution times are measured on the same MacBook described in the previous section.

An interesting metric is the break-even point, which is the
number of instances that need to be solved before
CVXPYgen is faster than CVXPY, when we include the code generation and compilation time.
This number is around 5000, and not too dependent on $N$.
A typical back-test might involve daily trading, with around 250 trading days in each year,
over 4 years, with hundreds of different hyper-parameter values, which gives
on the order of 100,000 solves, well above this break-even point.

\section{Conclusion}\label{sec:conclusion}

We have described CVXPYgen, a tool for generating custom C code that solves instances
of a family of convex optimization problems specified within CVXPY.
This gives a seemless path from prototyping an application using Python and CVXPY, 
to a final embedded implementation in C.
In addition to CVXPYgen supporting a wider variety of problems (such as SOCPs)
than the state-of-the-art code generator CVXGEN, numerical experiments
show that it outperforms CVXGEN in terms of
allowable problem size, compiled code size, and solve times.
For applications running on general purpose machines,
we obtain a significant speedup over CVXPY when many problem
instances are to be solved.

\appendices

\section{MPC example}\label{sec:appendix_MPC}

We use MPC to track the position and velocity of a quadcopter with mass $m$, experiencing gravitational acceleration $g$.
We model the quadcopter as a point mass with position error $p_k \in \reals^3$ and 
velocity error $v_k \in \reals^3$, where $k$ denotes the time step or period.
We concatenate the position and velocity error to the 
state $z_k = (p_k, v_k) \in \reals^6$, which we regulate to $0$.
The input is the force vector $u_k \in \reals^3$ without gravity compensation.
The dynamics are
\[
z_{k+1} = A z_k + B u_k,
\]
where $A \in \reals^{6 \times 6}$ and $B \in \reals^{6 \times 3}$.

We limit the tilt angle of the quadcopter.
Since the quadcopter's attitude is tied to the pointing direction of 
$u$ plus gravity compensation, we impose a (polyhedral) tilt angle constraint as
\[
c_j^T \left(u_k\right)_{0:1} \le \gamma \left(\left(u_k\right)_2  + mg\right),
\quad j = 0,\ldots,N^\text{hs}-1,
\]
where $\left(\cdot\right)_{0:1}$ and $\left(\cdot\right)_2$ denote the horizontal and 
vertical part of a vector in $\reals^3$ space, respectively, and $\gamma > 0$.
We use $N^\text{hs}$ halfspaces parametrized through $c_j$.
Compared to a spherical (natural) tilt angle constraint, when added to the MPC constraints, 
this formulation renders the problem QP representable.
The lower and upper thrust limits of the propellers are represented as 
$u_\text{vmin} \le \left(u_k\right)_2 \le u_\text{vmax}$ with $u_\text{vmin}<0$ and $u_\text{vmax} > 0$.

Up to horizon $H$, we penalize state errors and control effort via the traditional 
quadratic cost
\[
 z_H^T Q_T z_H + \sum_{k=0}^{H-1} \left( z_k^T Q z_k + u_k^T R u_k \right),
\]
with diagonal positive definite matrices $Q$ and $R$, and positive definite $Q_T$,
which is the solution to the discrete-time algebraic Riccati equation 
(as a function of $A$, $B$, $Q$, and $R$).
In addition, at every stage, we discourage rapid changes of the input 
(that the low-level attitude control system cannot follow) with the 
additional cost term
\[
\sum_{k=0}^{H-1}
\left(u_{k+1}-u_{k}\right)^T T \left(u_{k+1}-u_{k}\right),
\]
where $T$ is diagonal positive definite.
Combining all the above constraints and cost terms, we arrive at the MPC problem
\[
\begin{array}{ll}
	\text{minimize} & \displaystyle z_H^T Q_T z_H + \sum_{k=0}^{H-1} \big( z_k^T Q z_k + u_k^T R u_k + \\[-.75em]
	& \hspace{7em} (u_{k+1}-u_{k})^T T (u_{k+1}-u_{k}) \big) \\[.5em]
	\text{subject to} & z_0 = z_\text{meas}, \quad u_0 = u_\text{prev} \\
	& z_{k+1} = A z_k + B u_k,\quad k=0, \ldots, H-1 \\
	& u_\text{vmin} \le (u_k)_2 \le u_\text{vmax}, \quad k=1,\ldots, H-1 \\
	& c_j^T \left(u_k\right)_{0:1} \le \gamma \left(\left(u_k\right)_2  + mg\right),\\
& \quad j = 1,\ldots, N^\text{hs}-1, \quad k = 1, \ldots, H-1,
\end{array}
\]
where the states $z_k$ and inputs $u_k$ are optimization variables.
The current state measurement is $z_\text{meas}$ and the solution for the input at the first stage from the previous solve is $u_\text{prev}$.
In practice, these two would most certainly be the only parameters of the problem.
However, for demonstration purposes, we declare all other symbols (except for variables, all $c_j$, $N^\text{hs}$, and $H$) as parameters.
This problem formulation is not DPP-compliant, \eg, because of the multiplication of parameters $\gamma$, $m$, and $g$.

Before rewriting the problem in DPP-compliant form, we define the following convenience notation for matrix slicing. We use the zero-based counting scheme.
$M_{r:t,c:d}$ is the slice of some matrix $M$ from its $r$th to its $t$th row (included) and from its $c$th to its $d$th column (included). We slice full columns by omitting the row indices, \ie, $M_c$ is the $c$th column of $M$.
Finally, we write the DPP-compliant problem as
\[
\begin{array}{ll}
	\text{minimize} & \norm{Q_T^{1/2} Z_H}_2^2 + \norm{Q^{1/2} Z_{0:H-1}}_F^2 \\
	& + \norm{R^{1/2} U_{0:H-1}}_F^2 + \norm{T^{1/2} \left(U_{1:H} - U_{0:H-1}\right)}_F^2 \\[.5em]
	\text{subject to} & Z_0 = z_\text{meas}, \quad U_0 = u_\text{prev} \\
	& Z_{1:H} = A Z_{0:H-1} + B U_{0:H-1} \\
	& u_\text{vmin} \le U_{2,1:H-1} \le u_\text{vmax} \\
	& c_j^T U_{0:1,1:H-1} \le \gamma U_{2,1:H-1} + d,\\
& \quad j=0, \ldots, N^\text{hs}-1,
\end{array}
\]
where $Z \in \reals^{6 \times (H+1)}$ and $U \in \reals^{3 \times (H+1)}$ are the variables and contain the state and input vectors, respectively, for increasing stage count in their columns. The vector $d \in \reals^{H-1}$ contains $\gamma m g$ in all its entries. 
The problem is parametrized by 
\[
\begin{array}{ll}
Q_T^{1/2},~ Q^{1/2},~ R^{1/2},~ T^{1/2},~A, ~B,\\
 ~\gamma,~ d,~ u_\text{vmin},~
u_\text{vmax},~ z_\text{meas},~u_\text{prev}.
\end{array}
\]
In the expressions above, $M^{1/2}$ denotes any squareroot of the positive definite 
matrix $M$, \eg, the transposed Cholesky factor,
$\Vert \cdot \Vert_F$ denotes the Frobenius norm,
and the inequalities are elementwise.

\section{Portfolio optimization example}\label{sec:appendix_portfolio}

We search for a portfolio of holdings in $N$ assets and a cash balance.
The corresponding weights are $w \in \reals^{N+1}$, where the last entry 
represents the cash balance, and $\ones^T w = 1$.
We impose a leverage limit $\|w\|_1 \leq L$, where $L\geq 1$ is a parameter.

The return $r \in \reals^{N+1}$ has mean (or forecast) $\alpha \in \reals^{N+1}$,
so the expected portfolio return is $\alpha^T w$.
The risk or variance of the portfolio return is $w^T\Sigma w$, where
$\Sigma$ is the positive definite covariance matrix of the asset returns.

We consider two additional objective terms.
One is a trading or transaction cost 
$(\kappa^\text{tc})^T|w-w^\text{prev}|$, where the absolute value is 
elementwise, $w^\text{prev}$ is the previous period weights,
and $\kappa^\text{tc} \ge 0$ is a parameter.
The other is a short-selling cost $(\kappa^\text{sh})^T (w)_-$,
where $(w)_- = \max\{-w,0\}$ (elementwise), and $\kappa^\text{sh} \ge 0$
is a parameter.

The overall objective function is
\[
\alpha^T w
- \gamma^\text{risk} w^T \Sigma w
- \gamma^\text{tc} \left(\kappa^\text{tc}\right)^T |w-w^\text{prev}|
- \gamma^\text{sh} \left(\kappa^\text{sh}\right)^T (w)_-,
\]
where $\gamma^\text{risk}$,
$\gamma^\text{tc}$,
and $\gamma^\text{sh}$ are positive parameters that scale the 
risk, transaction cost, and shorting cost, respectively.

Our final optimization problem is
\[
	\begin{array}{ll}
		\text{maximize} \quad &\alpha^T w - \gamma^\text{risk} w^T \Sigma w \\
		& - \gamma^\text{tc} \left(\kappa^\text{tc}\right)^T |w-w^\text{prev}| - \gamma^\text{sh} \left(\kappa^\text{sh}\right)^T (w)_- \\[.5em]
		\text{subject to} \quad &\ones^T w = 1, \quad \norm{w}_1 \le L,
	\end{array}
\]
where $w \in \reals^{N+1}$ is the variable, and all other symbols are parameters.

The covariance matrix $\Sigma$ takes the standard factor model form,
\[
\Sigma = FF^T+D,
\]
where $F \in \reals^{(N+1) \times K}$ and $D$ is positive definite diagonal.
The number of factors in this risk model is $K$, which is usually much less than $N$.

This problem formulation is not DPP-compliant, but we can
rewrite it in DPP-compliant form by eliminating quadratic forms and collecting products 
of parameters, as in
\[
\begin{array}{ll}
	\text{maximize} & \big(\frac{\alpha}{\gamma^\text{risk}}\big)^T w - \norm{F^T w}_2^2 - \norm{D^{1/2} w}_2^2 \\
	& -\big(\frac{\gamma^\text{tc}}{\gamma^\text{risk}}\kappa^\text{tc}\big)^T |\Delta w| - \big(\frac{\gamma^\text{sh}}{\gamma^\text{risk}} \kappa^\text{sh}\big)^T (w)_- \\[.5em]
	\text{subject to} & \ones^T w = 1, \quad \norm{w}_1 \le L \\
	&\Delta w = w-w^\text{prev},
\end{array}
\]
where the portfolio weight vector $w \in \reals^{N+1}$ and the weight change vector $\Delta w \in \reals^{N+1}$ are the variables.
The problem is parametrized by
\[
\frac{\alpha}{\gamma^\text{risk}}, \quad
F, \quad
D^{1/2}, \quad
\frac{\gamma^\text{tc}}{\gamma^\text{risk}}\kappa^\text{tc}, \quad
\frac{\gamma^\text{sh}}{\gamma^\text{risk}}\kappa^\text{sh}, \quad
L, \quad
w^\text{prev}.
\]

We consider $N$ stock assets, chosen randomly from the S\&P 500,
with historical return data from 2017--2019.
For each value of $N$, we set $K=\max(N/10, 5)$.

\section*{Acknowledgment}

We would like to express our gratitude to John Lygeros for enabling this collaboration.
Moreover, we would like to thank JunEn Low and Mac Schwager for providing the quadcopter testing environment at the Stanford Flight Room.

\balance

\bibliography{refs}

\end{document}